\documentclass[11pt,amsfonts]{article}
\usepackage{graphicx}
\usepackage{latexsym}
\usepackage{amssymb}
\usepackage{amsmath}
\usepackage{stmaryrd}
\usepackage{layout}
\newtheorem{prop}{Proposition}
\newtheorem{lemma}{Lemma}
\newtheorem{definition}{Definition}

\newtheorem{theorem}{Theorem}
\newtheorem{remark}{Remark}

\def\real{{\mathord{{\rm I\kern-2.8pt R}}}}        % Fake blackboard bold R.
\def\inte{{\mathord{{\rm I\kern-2.8pt N}}}}

\def\sZZ{{\rm Z\kern-2.8ptem{}Z}}

\def\z{{\mathchoice
  {\sZZ}
  {\sZZ}
  {\rm Z\kern-0.30em{}Z}
  {\rm Z\kern-0.25em{}Z} }}
\def\sQQ{{\kern 0.27em \vrule height1.45ex width0.03em depth0em
          \kern-0.30em \rm Q}}
\def\qu{{\mathchoice
    {\sQQ}
    {\sQQ}
  {\kern 0.225em \vrule height1.05ex width0.025em depth0em \kern-0.25em \rm Q}
  {\kern 0.180em \vrule height0.78ex width0.020em depth0em \kern-0.20em \rm Q}
        }}
\def\sCC{{\kern 0.27em \vrule height1.45ex width0.03em depth0em
          \kern-0.30em \rm C}}
\def\complex{{\mathchoice
    {\sCC}
    {\sCC}
  {\kern 0.225em \vrule height1.05ex width0.025em depth0em \kern-0.25em \rm C}
  {\kern 0.180em \vrule height0.78ex width0.020em depth0em \kern-0.20em \rm C}
        }}

\newenvironment{dem}{{\bf Proof:}}{}

\newcommand{\ba}{\begin{array}}
\newcommand{\ea}{\end{array}}
\newcommand{\be}{\begin{equation}}
\newcommand{\ee}{\end{equation}}
\newcommand{\bea}{\begin{eqnarray}}
\newcommand{\eea}{\end{eqnarray}}
\newcommand{\beaa}{\begin{eqnarray*}}
\newcommand{\eeaa}{\end{eqnarray*}}

\newcommand{\eps}{\varepsilon}

\def\z{\zeta}

\font\tenmath=msbm10 \font\sevenmath=msbm7 \font\fivemath=msbm5
\newfam\mathfam \textfont\mathfam=\tenmath
\scriptfont\mathfam=\sevenmath \scriptscriptfont\mathfam=\fivemath

\def \={{\buildrel {\rm (law)} \over =}}

\def\qed{ \hfill \vrule width.25cm height.25cm depth0cm\smallskip}

\newcommand{\basa}{\begin{assumption}}
\newcommand{\easa}{\end{assumption}}

\newcommand{\bas}{\begin{assum}}
\newcommand{\eas}{\end{assum}}

\newcommand{\ignore}[1]{}
\begin{document}

\renewcommand{\thefootnote}{\fnsymbol{footnote}}

\renewcommand{\thefootnote}{\fnsymbol{footnote}}

\date{  }
\title{Berry-Ess\'een Bounds for Long Memory Moving Averages via Stein's Method and Malliavin Calculus}
\author{ Solesne Bourguin $^{1}\qquad $%
Ciprian A. Tudor $^{2,}$ \footnote{Associate member of the team Samm, Universit\'e de Paris 1 Panth\'eon-Sorbonne}\vspace*{0.1in} \\
$^{1}$ SAMM, Universit\'e de Paris 1 Panth\'eon-Sorbonne,\\
90, rue de Tolbiac, 75634, Paris, France. \\
solesne.bourguin@univ-paris1.fr \vspace*{0.1in} \\
 $^{2}$ Laboratoire Paul Painlev\'e, Universit\'e de Lille 1\\
 F-59655 Villeneuve d'Ascq, France.\\
 \quad tudor@math.univ-lille1.fr\vspace*{0.1in}}
\maketitle

\begin{abstract}
\noindent Using the Stein method on Wiener chaos introduced in \cite{NoPe1} we prove Berry-Ess\'een bounds for long memory moving averages.
\end{abstract}

\vskip0.3cm

{\bf 2010 AMS Classification Numbers:}  60F05, 60H05, 91G70.

 \vskip0.3cm


\section{Introduction}
Let $(\Omega, {\mathcal{F}}, P)$ be a probability space and $(W_{t}) _{t \geq 0}$ be a Brownian motion on this space. Let $F$ be a random variable defined on $\Omega$ which is differentiable in the sense of Malliavin calculus.
Then, using Stein's method on Wiener chaos, introduced by Nourdin and Peccati in \cite{NoPe1} (see also \cite{NoPe2} and \cite{NoPe3}), it is possible to measure the distance between the law of $F$ and the standard normal law $N(0,1)$. This distance can be defined in several ways (the Kolmogorov distance, the Wasserstein distance, the total variation distance or the Fortet-Mourier distance). More precisely we have, if ${\cal{L}}(F)$  denotes the law of $F$,
\begin{equation}\label{s1}
d({\cal{L}}(F), N(0,1)) \leq c\sqrt{ E\left( 1-\langle DF , D(-L) ^{-1} F \rangle _{L^{2}([0,1])}\right)^{2}}.
\end{equation}
Here, $D$ denotes the Malliavin derivative with respect to $W$ while $L$ is the generator of the Ornstein-Uhlenbeck semigroup. We will explain in the next section  how these operators are defined. The constant $c$ is equal to 1 in the case of the Kolmogorov and of the Wasserstein distance, $c$=2 for the  total variation distance and $c=4$ in the case of the Fortet-Mourier distance.
\\\\
\noindent These results have already been used to prove error bounds in various central limit theorems. In \cite{NoPe1}  the authors prove Berry-Ess\'een bounds in the central limit theorem for the subordinated functionals of the fractional Brownian motion and \cite{NoPe2} focuses on central limit theorems for Toeplitz quadratic functionals of continuous-time stationary processes. In \cite{NoPeRe} the authors extended the Stein's method to multidimensional settings. See also \cite{AMV}.
\\\\
\noindent In this paper  we will consider long memory moving averages defined by
\begin{equation*}
X_{n}= \sum _{i\geq 1} a_{i} \varepsilon _{n-i}, n\in \mathbb{Z}
\end{equation*}
where the innovations $\varepsilon _{i}$ are centered i.i.d. random variables having at least finite second moments and the moving averages $a_{i}$ are of the form $a_{i}=  i^{-\beta} L(i)$ with $\beta \in (\frac{1}{2},1)$ and $L$ slowly varying towards infinity. The covariance function $\rho(m)=\mathbf{E}\left(X_{0}X_{m}\right)$ behaves as $c_{\beta}m^{-2\beta +1}$ when $m\to \infty$ and consequently is not summable since $\beta < 1$. Therefore $X_{n}$ is usually called long-memory  or ``long-range dependence''  moving average. Let $K$ be a deterministic function which has Hermite rank $q$ and satisfies $\mathbf{E}(K^{2}(X_{n})) < \infty$ and define
\begin{equation*}
S_{N}=\sum_{n=1}^{N} \left[ K(X_{n})- \mathbf{E}\left(K(X_{n})\right)\right].
\end{equation*}
Then it has been proven in \cite{HH} (see also \cite{Wu}) that, with $c_{1}(\beta, q), c_{2}(\beta ,q)$ being positive constants depending only on $q$ and $\beta$: a) If $q>\frac{1}{2\beta -1}$, then the sequence $c_{1}(\beta, q)\frac{1}{\sqrt{N}} S_{N} $ converges in law to a standard normal random variable  and b) If $q<\frac{1}{2\beta -1} $, then the sequence $c_{2}(\beta, q)N^{\beta q - \frac{q}{2} - 1} S_{N}$ converges in law to a Hermite random variable of order $q$. This Hermite random variable, which will be defined in the next section,  is actually an iterated integral of a deterministic function with $q$ variables with respect to a Wiener process. This theorem is a variant of the non-central limit theorem from \cite{DM} and \cite{Ta2}.
In order to apply the techniques based on the Malliavin calculus and multiple Wiener-It\^o integrals, we will restrict our focus to the following situation: the innovations $\varepsilon _{i}$ are chosen to be the increments of a Brownian motion $W$ on the real line while the function $K$ is a Hermite polynomial of order $q$.  In this case the random variable $X_{n}$ is  a Wiener integral with respect to $W$, and $H_{q}(X_{n})$ can be expressed as a multiple Wiener-It\^o stochastic integral of order $q$ with respect to $W$. When $q>\frac{1}{2\beta -1}$ we will apply formula (\ref{s1}) in order to obtain the rate of convergence of $S_{N}$. When $q<\frac{1}{2\beta -1}$ the limit  of $S_{N}$ (after normalization) is not Gaussian and so we will use a different argument based on a result in \cite{DaMa} that has already been exploited in \cite{BN}.
\\\\
\noindent The paper is organized as follows. Section 2 deals with notation and preliminaries, such as the definition of a moving average process and a Wiener process on $\mathbb{R}$, but also gives a brief introduction to the tools of Malliavin calculus. In section 3, we will prove the Berry-Ess\'een bounds for the central and non central limit theorems for long-memory moving averages. Section 4 shows an application of our results to the Hsu-Robbins and Spitzer theorems for moving averages.

\section{Notation and Preliminaries}
In this section, we will give the main properties of infinite moving average processes and a proper definition of a Brownian motion on $\mathbb{R}$. We will relate one to the other to prove that the processes that we will consider in the latter are well defined. To conclude the preliminaries, we will finally focus on the sequences and results, such as central and non-central limit theorems that interest us in this paper.
\subsection{The Infinite Moving Average Process}
Before introducing the infinite moving average process, we will need the proper definition of a white noise on $\mathbb{Z}$.
\begin{definition}
\label{whitenoisedef}
The process $\left\{Z_{t}\right\}_{t \in \mathbb{Z}}$ is said to be a white noise with zero mean and variance $\sigma^{2}$, written $$\left\{Z_{t}\right\} \sim \mathcal{WN}(0,\sigma^{2}),$$ if and only if for every $h \in \mathbb{N}$, $\left\{Z_{t}\right\}$ has zero mean and covariance function $\gamma(h) = \mathbf{E}\left(Z_{t+h}Z_{t}\right)$ defined by $$
\gamma(h) =\left\{
\begin{array}{rl}
\sigma^{2} \mbox{\   \   if \   } h = 0 \\
0 \mbox{\   \   if \   } h \neq 0.
\end{array}
\right.
$$
\end{definition}
Now we can define the infinite moving average process.
\begin{definition}
If $\left\{Z_{t}\right\} \sim \mathcal{WN}(0,\sigma^{2})$ then we say that $\left\{X_{t}\right\}$ is a moving average $(\mbox{MA}(\infty))$ of $\left\{Z_{t}\right\}$ if there exists a sequence $\left\{\psi_{j}\right\}$ with $\sum_{j = 0}^{\infty}\left|\psi_{j}\right| < \infty$ such that
\begin{equation}
\label{movingav}
X_{t} = \sum_{j = 0}^{\infty}\psi_{j}Z_{t-j}, \mbox{\   \   \   } t = 0, \pm 1, \pm 2, ...
\end{equation}
\end{definition}
We have the following proposition on infinite moving averages (see \cite{BrDa} p. 91).
\begin{prop}
The $(\mbox{MA}(\infty))$ process defined by (\ref{movingav}) is stationary with mean zero and covariance function
\begin{equation}
\gamma(k) = \sigma^{2}\sum_{j = 0}^{\infty}\psi_{j}\psi_{j+\left|k\right|}.
\end{equation}
\end{prop}
For further details on moving averages, see \cite{BrDa} for a complete survey of this topic.
\subsection{The Brownian Motion on $\mathbb{R}$}
Here, we will give a proper definition of a two-sided Brownian motion on $\mathbb{R}$ (as defined in \cite{Ch}). We will then connect this definition to the underlying Hilbert space.
\begin{definition}
A two sided Brownian motion $\left\{W_{t}\right\}_{t \in \mathbb{R}}$ on $\mathbb{R}$ is a continous centered Gaussian process with covariance function
\begin{equation}
R(t,s) = \frac{1}{2}\left(\left|s\right| + \left|t\right| - \left|t-s\right|\right), \   \   s,t \in \mathbb{R}.
\end{equation}
Let $\mathcal{H} = L^{2}(\mathbb{R})$ be the underlying Hilbert space of this particular process. We have
\begin{equation}
R(t,s) = \left\{
\begin{array}{ll}
\left\langle \mathbf{1}_{\left[0,s\right]}, \mathbf{1}_{\left[0,t\right]}\right\rangle_{\mathcal{H}} = s \wedge t \mbox{\  if \  } s,t \geq 0 \\
\left\langle \mathbf{1}_{\left[s,0\right]}, \mathbf{1}_{\left[0,t\right]}\right\rangle_{\mathcal{H}} = 0 \mbox{\  if \  } s\leq 0 \mbox{\  and \  } t \geq 0 \\
\left\langle \mathbf{1}_{\left[0,s\right]}, \mathbf{1}_{\left[t,0\right]}\right\rangle_{\mathcal{H}} = 0 \mbox{\  if \  } s\geq 0 \mbox{\  and \  } t \leq 0 \\
\left\langle \mathbf{1}_{\left[s,0\right]}, \mathbf{1}_{\left[t,0\right]}\right\rangle_{\mathcal{H}} = - (s \vee t) = \left|s\right| \wedge \left|t\right| \mbox{\  if \  } s,t \leq 0.
\end{array}
\right.
\end{equation}
We could also define the two-sided Brownian motion by considering two independent standard Brownian motions on $\mathbb{R}^{+}$, $\left\{W_{t}^{(1)}\right\}$ and $\left\{W_{t}^{(2)}\right\}$ and by setting
\begin{equation}
W_{t} = \left\{
\begin{array}{ll}
W_{t}^{(1)} \mbox{\  if \  } t \geq 0 \\
W_{-t}^{(2)} \mbox{\  if \  } t \leq 0.
\end{array}
\right.
\end{equation}
$\left\{W_{t}\right\}$ has the same law as the one induced by the first definition.
\end{definition}
If we define the process $\left\{I_{t}\right\}_{t \in \mathbb{Z}}$ as the increment of the two-sided Brownian motion between $t$ and $t+1$, $t \in \mathbb{Z}$, we have $I_{t} = W_{t+1} - W_{t}$. The following holds.
\begin{prop}
The process $\left\{I_{t}\right\}_{t \in \mathbb{Z}}$ is a white noise on $\mathbb{Z}$ with mean 0 and variance 1.
\end{prop}
\begin{dem}
It is clear that $\left\{I_{t}\right\}$ is a centered Gaussian process. We only need to verify its covariance function. We have, for every $h \in \mathbb{Z}$,
\begin{equation*}
\mathbf{E}\left(I_{t+h}I_{t}\right) = \mathbf{E}\left((W_{t+h+1} - W_{t+h})(W_{t+1} - W_{t})\right) = \left\{
\begin{array}{ll}
1 \mbox{\  if \  } h = 0 \\
0 \mbox{\  if \  } h \neq 0
\end{array}
\right.
\end{equation*}
\qed
\end{dem}
\subsection{Limit Theorems for Functionals of i.i.d Gaussian Processes}
Here, we will focus on the following type of sequences
\begin{equation}\label{sn}
S_{N}= \sum_{n=1}^{N} \left[K(X_{n}) - \mathbf{E}\left(K(X_{n})\right)\right]
\end{equation}
where
\begin{equation}
X_{n}= \sum_{i=1}^{\infty} \alpha_{i}\left(W_{n-i} - W_{n-i-1}\right),
\end{equation}
with $\alpha_{i} \in \mathbb{R}$ and $\sum_{i=1}^{\infty} \alpha_{i}^{2} = 1$.
Note that $\left\{X_{n}\right\}$ is an infinite moving average of the white noise $\left\{I_{t}\right\} = \left\{W_{t+1} - W_{t}\right\}$. Thus its covariance function is given by
\begin{equation}
\label{covarXn}
\rho(m) := \sum_{i = 1}^{\infty}\alpha_{i}\alpha_{i+\left|m\right|}.
\end{equation}
For those sequences, central and non-central limit theorems have been proven. Here are the main results we will be focusing on.
\begin{theorem}
\label{CLTGen}
Suppose that the $\alpha_{i}$ are regularly varying with exponent $-\beta$, $\beta \in (1/2, 1)$ (i.e. $\alpha_{i} = \left|i\right|^{-\beta}L(i)$ and that $L(i)$ is slowly varying at $\infty$). Suppose that $K$ has Hermite rank $k$ and satisfies $\mathbf{E}(K^{2}(X_{n})) < \infty$. Then
\begin{description}
\item{i. }If $k<(2\beta - 1)^{-1}$, then
\begin{equation}
h_{k,\beta}^{-1}N^{\beta q-\frac{q}{2}-1}S_{N} \underset{N \rightarrow +\infty}{\longrightarrow} Z^{(k)}
\end{equation}
where $Z^{(k)}$ is a Hermite random variable of order $k$ defined by (\ref{hermite}) and $h_{k,\beta}$ is a positive constant depending on $k$ and $\beta$ (which will be defined later by (\ref{hq})).
\item{ii. }If $k>(2\beta - 1)^{-1}$, then
\begin{equation}
\frac{1}{\sigma _{k,\beta }\sqrt{N}}S_{N} \underset{N \rightarrow +\infty}{\longrightarrow} \mathcal{N}(0,1)
\end{equation}
with $\sigma _{k, \beta }$ defined by (\ref{sigma}).
\end{description}
\end{theorem}

\noindent We will compute the Berry-Ess\'een bounds for these central limit (CLT) and non-central limit (NCLT) theorems using Stein's Method and Malliavin Calculus. In the next paragraph, we will give the basic elements on these topics.

\subsection{Multiple Wiener-It\^o Integrals and Malliavin Derivatives}
Here we describe the elements from stochastic analysis that we will need in the paper. Consider ${\mathcal{H}}$ a real separable Hilbert space and $(B (\varphi), \varphi\in{\mathcal{H}})$ an isonormal Gaussian process on a probability space $(\Omega, {\cal{A}}, P)$, which is a centered Gaussian family of random variables such that $\mathbf{E}\left( B(\varphi) B(\psi) \right)  = \langle\varphi, \psi\rangle_{{\mathcal{H}}}$. Denote by $I_{n}$ the multiple stochastic integral with respect to
$B$ (see \cite{N}). This $I_{n}$ is actually an isometry between the Hilbert space ${\mathcal{H}}^{\odot n}$(symmetric tensor product) equipped with the scaled norm $\frac{1}{\sqrt{n!}}\Vert\cdot\Vert_{{\mathcal{H}}^{\otimes n}}$ and the Wiener chaos of order $n$ which is defined as the closed linear span of the random variables $H_{n}(B(\varphi))$ where $\varphi\in{\mathcal{H}}, \Vert\varphi\Vert_{{\mathcal{H}}}=1$ and $H_{n}$ is the Hermite polynomial of degree $n\geq 1$
\begin{equation*}
H_{n}(x)=\frac{(-1)^{n}}{n!} \exp \left( \frac{x^{2}}{2} \right)
\frac{d^{n}}{dx^{n}}\left( \exp \left( -\frac{x^{2}}{2}\right)
\right), \hskip0.5cm x\in \mathbb{R}.
\end{equation*}
The isometry of multiple integrals can be written as: for $m,n$ positive integers,
\begin{eqnarray}
\mathbf{E}\left(I_{n}(f) I_{m}(g) \right) &=& n! \langle f,g\rangle _{{\mathcal{H}}^{\otimes n}}\quad \mbox{if } m=n,\nonumber \\
\mathbf{E}\left(I_{n}(f) I_{m}(g) \right) &= & 0\quad \mbox{if } m\not=n.\label{iso}
\end{eqnarray}
It also holds that
\begin{equation*}
I_{n}(f) = I_{n}\big( \tilde{f}\big)
\end{equation*}
where $\tilde{f} $ denotes the symmetrization of $f$ defined by $\tilde{f}%
(x_{1}, \ldots , x_{n}) =\frac{1}{n!} \sum_{\sigma \in {\cal S}_{n}}
f(x_{\sigma (1) }, \ldots , x_{\sigma (n) } ) $.
\\\\
We recall that any square integrable random variable which is measurable with respect to the $\sigma$-algebra generated by $B$ can be expanded into an orthogonal sum of multiple stochastic integrals
\begin{equation}
\label{sum1} F=\sum_{n\geq0}I_{n}(f_{n})
\end{equation}
where $f_{n}\in{\mathcal{H}}^{\odot n}$ are (uniquely determined)
symmetric functions and $I_{0}(f_{0})=\mathbf{E}\left[  F\right]$.
\\\\
Let $L$ be the Ornstein-Uhlenbeck operator
\begin{equation*}
LF=-\sum_{n\geq 0} nI_{n}(f_{n})
\end{equation*}
if $F$ is given by (\ref{sum1}).
\\\\
For $p>1$ and $\alpha \in \mathbb{R}$ we introduce the Sobolev-Watanabe space $\mathbb{D}^{\alpha ,p }$  as the closure of
the set of polynomial random variables with respect to the norm
\begin{equation*}
\Vert F\Vert _{\alpha , p} =\Vert ((I -L)F) ^{\frac{\alpha }{2}} \Vert_{L^{p} (\Omega )}
\end{equation*}
where $I$ represents the identity. We denote by $D$  the Malliavin  derivative operator that acts on smooth functions of the form $F=g(B(\varphi _{1}), \ldots , B(\varphi_{n}))$ ($g$ is a smooth function with compact support and $\varphi_{i} \in {{\cal{H}}}$)
\begin{equation*}
DF=\sum_{i=1}^{n}\frac{\partial g}{\partial x_{i}}(B(\varphi _{1}), \ldots , B(\varphi_{n}))\varphi_{i}.
\end{equation*}
The operator $D$ is continuous from $\mathbb{D}^{\alpha , p} $ into $\mathbb{D} ^{\alpha -1, p} \left( {\cal{H}}\right).$
\\\\
\noindent In this paper we will use the Malliavin calculus with respect to the Brownian motion on $\mathbb{R}$ as introduced above. Note that the Brownian motion on the real line is an isonormal process and its underlying Hilbert space is ${\cal{H}}=L^{2}(\mathbb{R})$.
\\\\
\noindent We will now introduce the Hermite random variable, which is the limit in Theorem 1, point i. The Hermite random variable of order $q $ is given by
\begin{equation}\label{hermite}
Z^{(q)} =d(q, \beta ) I_{q} (g(\cdot ) )
\end{equation}
where
\begin{equation}\label{kernel}g(y_{1},.., y_{q})= \int_{y_{1}\vee ... \vee y_{q}}^{1}du (u-y_{1})_{+}^{-\beta } ...(u-y_{q})_{+}^{-\beta }.
\end{equation}
The constant $d(q, \beta)$ is a normalizing constant which ensures that $\mathbf{E}(Z^{(q)})^{2} =1$. This constant is explicitly computed below.
\begin{eqnarray*}
\mathbf{E}(Z^{(q)})^{2}&=& q! d(q, \beta)^{2} \int_{0}^{1} dudv \left( \int_{\mathbb{R}} (u-y)_{+}^{-\beta} (v-y)_{+}^{-\beta} dy \right) ^{q} \\
&=&q! d(q, \beta)^{2} \beta (2\beta-1, 1-\beta) ^{q} \int_{0}^{1} dudv \vert u-v\vert ^{-2qb+q}\\
&=&q! d(q, \beta)^{2} \beta (2\beta-1, 1-\beta) ^{q} \frac{2}{(-2\beta q+q+1)(-2\beta q+q+2)}
\end{eqnarray*}
where we used
$$\int_{\mathbb{R}} (u-y)_{+}^{-\beta} (v-y)_{+}^{-\beta} dy=\beta (2\beta-1, 1-\beta)\vert u-v\vert ^{-2\beta +1}=c_{\beta}\vert u-v\vert ^{-2\beta +1}$$
and we denoted $c_{\beta}:=\boldsymbol{\beta}(2\beta-1, 1-\beta)$, $\boldsymbol{\beta}$ being the beta function defined by $$\boldsymbol{\beta}(x,y) = \int_{0}^{1}t^{x-1}(1-t)^{y-1}dt = \int_{0}^{\infty}\frac{t^{x-1}}{(1+t)^{x+y}}dt.$$
Therefore
\begin{equation}\label{dq}
d(q, \beta)^{2}=\frac{(-2\beta q+q+1)(-2\beta q+q+2)}{2q! c_{\beta} ^{q}}.
\end{equation}

\subsection{Stein's Method on a Fixed Wiener Chaos}
Let $F = I_{q}(h), h \in \mathcal{H}^{\odot q}$ be an element on the Wiener chaos of order $q$. Recall that for any fixed $z \in \mathbb{R}$, the Stein equation is given by
\begin{eqnarray}
\label{steinEq}
\mathbf{1}_{\left(_\infty,z\right]}(x) - \Phi(x) = f'(x) - xf(x).
\end{eqnarray}
It is well known that (\ref{steinEq}) admits a solution $f_{z}$ bounded by $\sqrt{2\pi}/4$ and such that $\left\|f'_{z}\right\|_{\infty} \leq 1$. By taking $x = F$ in (\ref{steinEq}) and by taking the expectation, we get
\begin{eqnarray}
\label{steinEqEsp}
\mathbf{P}(F \leq z) - \mathbf{P}(N \leq z) = \mathbf{E}\left(f_{z}'(F) - Ff_{z}(F)\right)
\end{eqnarray}
where $N$ is a standard normal random variable ($N \hookrightarrow \mathcal{N}(0,1)$). By writing $F = LL^{-1}F = -\delta DL^{-1}F$ and by integrating by part, we find
\begin{eqnarray}
\mathbf{E}\left(Ff_{z}(F)\right) &=& \mathbf{E}\left(-\delta DL^{-1}Ff_{z}(F)\right) = \mathbf{E}\left(\left\langle - DL^{-1}F,D(f_{z}(F))\right\rangle_{\mathcal{H}}\right)\nonumber
\\
&=& \mathbf{E}\left(\left\langle - DL^{-1}F,f_{z}'(F)DF\right\rangle_{\mathcal{H}}\right) = \mathbf{E}\left(f_{z}'(F)\left\langle - DL^{-1}F,DF\right\rangle_{\mathcal{H}}\right).\nonumber
\end{eqnarray}
Thus, by replacing in (\ref{steinEqEsp}), we obtain
\begin{eqnarray}
\mathbf{E}\left(f_{z}'(F) - Ff_{z}(F)\right) = \mathbf{E}\left(f_{z}'(F)\left(1 - \left\langle - DL^{-1}F,DF\right\rangle_{\mathcal{H}}\right)\right)\nonumber
\end{eqnarray}
and
\begin{eqnarray}
\mathbf{P}(F \leq z) - \mathbf{P}(N \leq z) = \mathbf{E}\left(f_{z}'(F)\left(1 - \left\langle - DL^{-1}F,DF\right\rangle_{\mathcal{H}}\right)\right).
\end{eqnarray}
On a different but related matter, the Kolmogorov distance is defined by
\begin{eqnarray}
d_{\mbox{\tiny{Kol}}}(X,Y) = \underset{z \in \mathbb{R}}{\mbox{sup}} \left|\mathbf{P}(X \leq z) - \mathbf{P}(Y \leq z)\right|.
\end{eqnarray}
Therefore we have
\begin{eqnarray}
d_{\mbox{\tiny{Kol}}}(F,N) = \underset{z \in \mathbb{R}}{\mbox{sup}} \left|\mathbf{E}\left(f_{z}'(F)\left(1 - \left\langle - DL^{-1}F,DF\right\rangle_{\mathcal{H}}\right)\right)\right|. \nonumber
\end{eqnarray}
By applying the Cauchy-Schwarz inequality, we get
\begin{eqnarray}
\label{dkol}
d_{\mbox{\tiny{Kol}}}(F,N) & \leq & \underbrace{\left[\mathbf{E}\left((f_{z}'(F))^{2}\right)\right]^{\frac{1}{2}}}_{\leq 1}\left[\mathbf{E}\left(\left(1 - \left\langle - DL^{-1}F,DF\right\rangle_{\mathcal{H}}\right)^{2}\right)\right]^{\frac{1}{2}} \nonumber
\\
& \leq & \sqrt{\mathbf{E}\left(\left(1 - \left\langle - DL^{-1}F,DF\right\rangle_{\mathcal{H}}\right)^{2}\right)}.
\end{eqnarray}
Recall that $F = I_{q}(h)$ and so in that case the equality $$\left\langle - DL^{-1}F,DF\right\rangle_{\mathcal{H}} = q^{-1}\left\|DF\right\|_{\mathcal{H}}^{2}$$ stands. Thus, we can rewrite (\ref{dkol}) as
\begin{eqnarray}
\label{dkolHq}
d_{\mbox{\tiny{Kol}}}(F,N) & \leq & c\sqrt{\mathbf{E}\left(\left(1 - q^{-1}\left\|DF\right\|_{\mathcal{H}}^{2}\right)^{2}\right)}
\end{eqnarray}
with $c=1$. As we mentioned in the introduction, the above inequality still holds true for other distances (Wasserstein, total variation or Fortet-Mourier). The constant $c$ is equal to 1 in the case of the Kolmogorov and of the Wasserstein distance, $c$=2 for the  total variation distance and $c=4$ in the case of the Fortet-Mourier distance.

\section{Berry-Ess\'een Bounds in the Central and Non-Central Limit Theorems}
As previously mentionned in the introduction, we will focus on the case where $K = H_{q}$, $H_{q}$ being the Hermite polynomial of order $q$. In this case, we will be able to give a more appropriate representation of $S_{N}$ in terms of multiple stochastic integrals. We will also assume $a_{i}=i^{-\beta}$ for every $i\geq 1$ so the slowly varying function at $\infty$ is chosen to be identically equal to one.

\subsection{Representation of $S_{N}$ as an Element of the $q^{\mbox{\tiny{th}}}$-Chaos}
Note that $X_{n}$ can also be written as
\begin{eqnarray}
X_{n} & = & \sum_{i=1}^{\infty} \alpha_{i}\left(W_{n-i} - W_{n-i-1}\right) = \sum_{i=1}^{\infty} \alpha_{i}I_{1}\left(\mathbf{1}_{\left[n-i-1,n-i\right]}\right) \nonumber
\\
& = & I_{1}\left(\underbrace{\sum_{i=1}^{\infty} \alpha_{i}\mathbf{1}_{\left[n-i-1,n-i\right]}}_{f_{n}}\right) =I_{1}\left(f_{n}\right).
\end{eqnarray}
As $K = H_{q}$, we have
\begin{eqnarray}
S_{N} & = & \sum_{n=1}^{N} \left[H_{q}(X_{n}) - \mathbf{E}\left(H_{q}(X_{n})\right)\right] =\sum_{n=1}^{N} \left[H_{q}(I_{1}(f_{n})) - \mathbf{E}\left(H_{q}(I_{1}(f_{n}))\right)\right] \nonumber
\end{eqnarray}
We know that, if $\left\|f\right\|_{\mathcal{H}} = 1$, we have $H_{q}(I_{1}(f)) = \frac{1}{q!}I_{q}(f^{\otimes q})$. Furthermore, we have
\begin{eqnarray}
\left\|f_{n}\right\|_{\mathcal{H}}^{2} & = & \left\langle f_{n},f_{n}\right\rangle_{\mathcal{H}} = \left\langle \sum_{i=1}^{\infty} \alpha_{i}\mathbf{1}_{\left[n-i-1,n-i\right]},\sum_{r=1}^{\infty} \alpha_{r}\mathbf{1}_{\left[n-r-1,n-r\right]}\right\rangle_{\mathcal{H}} \nonumber
\\
& = & \sum_{i,r=1}^{\infty}\alpha_{i}\alpha_{r} \left\langle \mathbf{1}_{\left[n-i-1,n-i\right]},\mathbf{1}_{\left[n-r-1,n-r\right]}\right\rangle_{\mathcal{H}}.\nonumber
\end{eqnarray}
It is easily verified that if $i>r \Leftrightarrow n-i \leq n-r-1$ or $i<r \Leftrightarrow n-r \leq n-i-1$, we have $\left[n-i-1,n-i\right]\cap \left[n-r-1,n-r\right] = \emptyset$ and thus $\left\langle \mathbf{1}_{\left[n-i-1,n-i\right]},\mathbf{1}_{\left[n-r-1,n-r\right]}\right\rangle_{\mathcal{H}} = 0$. It follows that
\begin{eqnarray}
\left\|f_{n}\right\|_{\mathcal{H}}^{2} & = & \sum_{i=1}^{\infty}\alpha_{i}^{2} \left\| \mathbf{1}_{\left[n-i-1,n-i\right]}\right\|_{\mathcal{H}}^{2} = \sum_{i=1}^{\infty}\alpha_{i}^{2} = 1.\nonumber
\end{eqnarray}
Thanks to this result, $S_{N}$ can be represented as
\begin{eqnarray}
S_{N} & = & \sum_{n=1}^{N} \left[H_{q}(I_{1}(f_{n})) - \mathbf{E}\left(H_{q}(I_{1}(f_{n}))\right)\right] = \frac{1}{q!}\sum_{n=1}^{N} \left[I_{q}(f_{n}^{\otimes q}) - \mathbf{E}\left(I_{q}(f_{n}^{\otimes q})\right)\right] \nonumber
\\
& = & \frac{1}{q!}\sum_{n=1}^{N} I_{q}(f_{n}^{\otimes q}) = \frac{1}{q!}I_{q}(\sum_{n=1}^{N}f_{n}^{\otimes q}).\nonumber
\end{eqnarray}

\subsection{Berry-Ess\'een Bounds for the Central Limit Theorem}
We will first focus on the case where $q > (2\beta - 1)^{-1}$, i.e. the central limit theorem. Let $Z_{N} = \frac{1}{\sigma\sqrt{N}}S_{N}$ where $\sigma _{q, \beta}$ is given by
\begin{eqnarray}\label{sigma}
\sigma:=\sigma_{q,\beta} ^{2} = \frac{1}{q!}\sum_{m = -\infty}^{+\infty}\left(\sum_{i = 1}^{\infty}\alpha_{i}\alpha_{i+\left|m\right|}\right)^{q} = \frac{1}{q!}\sum_{m = -\infty}^{+\infty}\rho^{q}(m).
\end{eqnarray}
The following result gives the Berry-Ess\'een bounds for the central limit part of theorem \ref{CLTGen}.
\begin{theorem}\label{tt1}
Under the condition $q > (2\beta -1)^{-1}$, $Z_{N}$ converges in law towards $Z \sim \mathcal{N}(0,1)$. Moreover, there exists a constant $C_{\beta}$, depending uniquely on $\beta$, such that, for any $N \geq 1$,
\begin{eqnarray*}
\underset{z \in \mathbb{R}}{\mbox{sup}} \left|\mathbf{P}(Z_{N} \leq z) - \mathbf{P}(Z \leq z)\right| \leq C_{\beta}\left\{
\begin{array}{rl}
N^{\frac{q}{2}+ \frac{1}{2} -q\beta} \mbox{\   \   if \   } \beta \in \left(\frac{1}{2},\frac{q}{2q-2}\right] \\
N^{\frac{1}{2}-\beta} \mbox{\   \   if \   } \beta \in \left[\frac{q}{2q-2},1\right)
\end{array}
\right.
\end{eqnarray*}
\end{theorem}
{\bf Proof: }
Theorem \ref{CLTGen} states that $Z_{N} \underset{N \rightarrow +\infty}{\longrightarrow} \mathcal{N}(0,1)$. Because of (\ref{dkol}) and (\ref{dkolHq}), we will evaluate the quantity $$\mathbf{E}\left(\left(1 - q^{-1}\left\|DZ_{N}\right\|_{\mathcal{H}}^{2}\right)^{2}\right).$$ We will start by computing $\left\|DZ_{N}\right\|_{\mathcal{H}}^{2}$. We have the following lemma.
\begin{lemma}
\label{l1new}
The following result on $\left\|DZ_{N}\right\|_{\mathcal{H}}$ holds.
\begin{eqnarray}
\frac{1}{q}\left\|DZ_{N}\right\|_{\mathcal{H}}^{2} - 1 &=& \sum_{r = 0}^{q-1}A_{r}(N) - 1 \nonumber
\end{eqnarray}
where
\begin{eqnarray}
A_{r}(N) = \frac{q r!}{\sigma^{2} N}\left(\begin{array}{cc}q-1 \\r \end{array}\right)^{2}\sum_{k,l = 1}^{N}I_{2q-2 - 2r}\left(f_{k}^{\otimes q-1-r} \widetilde{\otimes}  f_{l}^{\otimes q-1-r}\right)\left\langle f_{k},f_{l} \right\rangle_{\mathcal{H}}^{r+1}.
\end{eqnarray}
\end{lemma}
{\bf Proof: }
We have
\begin{eqnarray*}
DZ_{N} &=& D\left(\frac{1}{\sigma\sqrt{N}}\sum_{n = 1}^{N}I_{q}\left(f_{n}^{\otimes q}\right)\right)= \frac{q}{\sigma\sqrt{N}}\sum_{n = 1}^{N}I_{q-1}\left(f_{n}^{\otimes q-1}\right)f_{n}
\end{eqnarray*}
and
\begin{eqnarray}
\label{normdzn}
\left\|DZ_{N}\right\|_{\mathcal{H}}^{2} &=& \frac{q^{2}}{\sigma^{2} N}\sum_{k,l = 1}^{N}I_{q-1}\left(f_{k}^{\otimes q-1}\right)I_{q-1}\left(f_{l}^{\otimes q-1}\right)\left\langle f_{k},f_{l} \right\rangle_{\mathcal{H}}.
\end{eqnarray}
The multiplication formula between multiple stochastic integrals gives us that
\begin{eqnarray}
I_{q-1}\left(f_{k}^{\otimes q-1}\right)I_{q-1}\left(f_{l}^{\otimes q-1}\right) &=& \sum_{r = 0}^{q-1}r!\left(\begin{array}{cc}q-1 \\r \end{array}\right)^{2}I_{2q-2 - 2r}\left(f_{k}^{\otimes q-1-r} \widetilde{\otimes}  f_{l}^{\otimes q-1-r}\right)\left\langle f_{k},f_{l} \right\rangle_{\mathcal{H}}^{r}. \nonumber
\end{eqnarray}
By replacing in (\ref{normdzn}), we obtain
\begin{eqnarray}
\left\|DZ_{N}\right\|_{\mathcal{H}}^{2} &=& \frac{q^{2}}{\sigma^{2} N}\sum_{r = 0}^{q-1}r!\left(\begin{array}{cc}q-1 \\r \end{array}\right)^{2}\sum_{k,l = 1}^{N}I_{2q-2 - 2r}\left(f_{k}^{\otimes q-1-r} \widetilde{\otimes}  f_{l}^{\otimes q-1-r}\right)\left\langle f_{k},f_{l} \right\rangle_{\mathcal{H}}^{r+1} \nonumber
\end{eqnarray}
and the conclusion follows easily.
\qed
\\\\
\noindent By using lemma \ref{l1new}, we can now evaluate $\mathbf{E}\left(\left(1 - q^{-1}\left\|DZ_{N}\right\|_{\mathcal{H}}^{2}\right)^{2}\right)$. We have
\begin{eqnarray}
\mathbf{E}\left(\left(1 - q^{-1}\left\|DZ_{N}\right\|_{\mathcal{H}}^{2}\right)^{2}\right) & = & \mathbf{E}\left(\left(\sum_{r = 0}^{q-1}A_{r}(N) - 1\right)^{2}\right) \nonumber
\\
& = & \mathbf{E}\left(\sum_{r = 0}^{q-1}A_{r}^{2}(N) + \sum_{r,p = 0, r \neq p}^{q-1}A_{r}(N)A_{p}(N) - 2\sum_{r = 0}^{q-1}A_{r}(N) + 1\right)\nonumber
\\
& = & \mathbf{E}\left(\sum_{r = 0}^{q-2}A_{r}^{2}(N) + A_{q-1}^{2}(N) + \sum_{r = 0}^{q-2}A_{r}(N)A_{q-1}(N)\right. \nonumber
\\
&& \left. + \sum_{r = 0}^{q-1}\sum_{p = 0, r \neq p}^{q-2}A_{r}(N)A_{p}(N) - 2\sum_{r = 0}^{q-1}A_{r}(N) + 1\right).\nonumber
\end{eqnarray}
Recall that $\mathbf{E}\left(I_{m}I_{n}\right) = 0$ if $m \neq n$. Thus, $\mathbf{E}\left(\sum_{r = 0}^{q-1}\sum_{p = 0, r \neq p}^{q-2}A_{r}(N)A_{p}(N)\right) = 0$. We can simplify our previous equality by writing
\begin{eqnarray}
\label{toev}
\mathbf{E}\left(\left(1 - q^{-1}\left\|DZ_{N}\right\|_{\mathcal{H}}^{2}\right)^{2}\right) & = & \mathbf{E}\left(\sum_{r = 0}^{q-2}A_{r}^{2}(N) + \sum_{r = 0}^{q-1}A_{r}(N)A_{q-1}(N)\right. \nonumber
\\
&& \left.  - \sum_{r = 0}^{q-1}A_{r}(N) - \sum_{r = 0}^{q-1}A_{r}(N) + 1\right)\nonumber
\\
& = & \mathbf{E}\left(\sum_{r = 0}^{q-2}A_{r}^{2}(N) + (A_{q-1}(N) - 1)\sum_{r = 0}^{q-1}A_{r}(N)\right. \nonumber
\\
&& \left. - \sum_{r = 0}^{q-2}A_{r}(N) - A_{q-1}(N) + 1\right)\nonumber
\\
& = & \mathbf{E}\left(\sum_{r = 0}^{q-2}A_{r}^{2}(N) - (A_{q-1}(N) - 1)\right)\nonumber
\\
& = & \sum_{r = 0}^{q-2}\mathbf{E}\left(A_{r}^{2}(N)\right) - (A_{q-1}(N) - 1).
\end{eqnarray}
We need to evaluate the behaviour of those two terms as $N \to \infty$, but first, recall that the $\alpha_{i}$ are of the form $\alpha_{i} = i^{-\beta}$ with $\beta \in \left(1/2, 1\right)$. We will use the notation $a_{n}\sim b_{n}$ meaning that $a_{n}$ and $b_{n}$ have the same limit as $n\to \infty$ and $a_{n} \trianglelefteqslant  b_{n}$ meaning that $\mbox{sup}_{n \geq 1} \left|a_{n}\right| / \left|b_{n}\right| < \infty$. Below is a useful lemma we will use throughout the paper.
\begin{lemma}\label{l1}
\begin{enumerate}
   \item We have $$\rho(n) \sim c_{\beta} n^{-2\beta + 1}$$
   with $c_{\beta}= \int_{0}^{\infty} y^{-\beta} (y+1)^{-\beta }dy = \boldsymbol{\beta}(2\beta-1,1-\beta).$ The constant $c_{\beta}$ is the same as the one in the definition of the Hermite random variable (see (\ref{dq})).
   \item For any $\alpha \in \mathbb{R}$, we have $$\sum_{k = 1}^{n-1}k^{\alpha} \trianglelefteqslant 1 + n^{\alpha + 1}.$$
   \item If $\alpha \in \left(-\infty , -1\right)$, we have $$\sum_{k = n}^{\infty}k^{\alpha} \trianglelefteqslant n^{\alpha + 1}.$$
\end{enumerate}
\end{lemma}
\begin{dem}
Points 2. and 3. follow from \cite{NoPe1}, Lemma 4.3. We will only prove the first point of the lemma (as the other points have been proven in \cite{NoPe1}). We know that $\rho(n) = \sum_{i = 1}^{\infty}i^{-\beta}\left(i+\left|n\right|\right)^{-\beta}$ behaves as $\int_{0}^{\infty}x^{-\beta}\left(x+\left|n\right|\right)^{-\beta}dx$ and the following holds
\begin{eqnarray*}
\int_{0}^{\infty}x^{-\beta}\left(x+\left|n\right|\right)^{-\beta}dx = \int_{0}^{\infty}x^{-\beta}\left|n\right|^{-\beta}\left(\frac{x}{\left|n\right|}+1\right)^{-\beta}dx = \left|n\right|^{-2\beta+1}\underbrace{\int_{0}^{\infty}y^{-\beta}\left(y+1\right)^{-\beta}dy}_{c_{\beta}}.
\end{eqnarray*}
Thus,
\begin{eqnarray*}
\rho(n) = \sum_{i = 1}^{\infty}i^{-\beta}\left(i+\left|n\right|\right)^{-\beta} \sim c_{\beta}n^{-2\beta+1}.
\end{eqnarray*}
\qed
\end{dem}
\\\\
\noindent We will start the evaluation of (\ref{toev}) with the term $(A_{q-1}(N) - 1)$ and we can write
\begin{eqnarray}
A_{q-1}(N) - 1 = \frac{q!}{\sigma^{2} N}\sum_{k,l = 1}^{N}\left\langle f_{k},f_{l} \right\rangle_{\mathcal{H}}^{q} - 1. \nonumber
\end{eqnarray}
Note that we have
\begin{eqnarray}
\left\langle f_{k},f_{l} \right\rangle_{\mathcal{H}} = \sum_{i = 1}^{\infty}\alpha_{i}\alpha_{i+\left|l-k\right|} = \rho(l-k). \nonumber
\end{eqnarray}
Hence
\begin{eqnarray}
\label{aqmoins1}
A_{q-1}(N) - 1 &=& \frac{q!}{\sigma^{2} N}\sum_{k,l = 1}^{N}\left(\sum_{i = 1}^{\infty}\alpha_{i}\alpha_{i+\left|l-k\right|}\right)^{q} - 1 \nonumber
\\
&=& \frac{1}{\sigma^{2} N}\left(q!\sum_{k,l = 1}^{N}\left(\sum_{i = 1}^{\infty}\alpha_{i}\alpha_{i+\left|l-k\right|}\right)^{q} - N\sigma^{2}\right) \nonumber
\\
&=& \frac{1}{\sigma^{2} N}\left(q!\sum_{k,l = 1}^{N}\left(\sum_{i = 1}^{\infty}\alpha_{i}\alpha_{i+\left|l-k\right|}\right)^{q} - Nq!\sum_{m = -\infty}^{+\infty}\left(\sum_{i = 1}^{\infty}\alpha_{i}\alpha_{i+\left|m\right|}\right)^{q}\right).
\end{eqnarray}
Observe that
\begin{eqnarray*}
\sum_{k,l = 1}^{N}\left(\sum_{i = 1}^{\infty}\alpha_{i}\alpha_{i+\left|l-k\right|}\right)^{q} &=& \sum_{k\leq l }^{N}\left(\sum_{i = 1}^{\infty}\alpha_{i}\alpha_{i+\left|l-k\right|}\right)^{q} + \sum_{k>l}^{N}\left(\sum_{i = 1}^{\infty}\alpha_{i}\alpha_{i+\left|l-k\right|}\right)^{q}
\\
&=& \sum_{k = 1 }^{N}\sum_{l = k }^{N}\left(\sum_{i = 1}^{\infty}\alpha_{i}\alpha_{i+\left|l-k\right|}\right)^{q} + \sum_{l = 1 }^{N}\sum_{k = l+1 }^{N}\left(\sum_{i = 1}^{\infty}\alpha_{i}\alpha_{i+\left|l-k\right|}\right)^{q}.
\end{eqnarray*}
Let $m = l-k$. We obtain
\begin{eqnarray*}
\sum_{k,l = 1}^{N}\left(\sum_{i = 1}^{\infty}\alpha_{i}\alpha_{i+\left|l-k\right|}\right)^{q} &=& \sum_{k = 1 }^{N}\sum_{m = 0 }^{N-k}\left(\sum_{i = 1}^{\infty}\alpha_{i}\alpha_{i+\left|m\right|}\right)^{q} + \sum_{l = 1 }^{N}\sum_{m = -N+l }^{-1}\left(\sum_{i = 1}^{\infty}\alpha_{i}\alpha_{i+\left|m\right|}\right)^{q}
\\
&=& \sum_{m = 0 }^{N-1}\sum_{k = 1 }^{N-m}\left(\sum_{i = 1}^{\infty}\alpha_{i}\alpha_{i+\left|m\right|}\right)^{q} + \sum_{m = -(N-1) }^{-1}\sum_{l = 1 }^{N+m}\left(\sum_{i = 1}^{\infty}\alpha_{i}\alpha_{i+\left|m\right|}\right)^{q}
\\
&=& \sum_{m = 0 }^{N-1}(N-m)\left(\sum_{i = 1}^{\infty}\alpha_{i}\alpha_{i+\left|m\right|}\right)^{q} + \sum_{m = -(N-1) }^{-1}(N+m)\left(\sum_{i = 1}^{\infty}\alpha_{i}\alpha_{i+\left|m\right|}\right)^{q}
\\
&=& N\sum_{m = -(N-1) }^{N-1}\left(\sum_{i = 1}^{\infty}\alpha_{i}\alpha_{i+\left|m\right|}\right)^{q} - 2\sum_{m =0 }^{N-1}m\left(\sum_{i = 1}^{\infty}\alpha_{i}\alpha_{i+\left|m\right|}\right)^{q}.
\end{eqnarray*}
By replacing in (\ref{aqmoins1}), we get
\begin{eqnarray*}
A_{q-1}(N) - 1 &=& \frac{q!}{\sigma^{2} N}\left(N\sum_{m = -(N-1) }^{N-1}\left(\sum_{i = 1}^{\infty}\alpha_{i}\alpha_{i+\left|m\right|}\right)^{q} - N\sum_{m = -\infty}^{+\infty}\left(\sum_{i = 1}^{\infty}\alpha_{i}\alpha_{i+\left|m\right|}\right)^{q} \right.
\\
&& \left. - 2\sum_{m =0 }^{N-1}m\left(\sum_{i = 1}^{\infty}\alpha_{i}\alpha_{i+\left|m\right|}\right)^{q}\right)
\\
&=& \frac{q!}{\sigma^{2} N}\left(-N\sum_{m = -\infty }^{-N}\left(\sum_{i = 1}^{\infty}\alpha_{i}\alpha_{i+\left|m\right|}\right)^{q} -N\sum_{m = N }^{\infty}\left(\sum_{i = 1}^{\infty}\alpha_{i}\alpha_{i+\left|m\right|}\right)^{q} \right.
\\
&& \left. - 2\sum_{m =0 }^{N-1}m\left(\sum_{i = 1}^{\infty}\alpha_{i}\alpha_{i+\left|m\right|}\right)^{q} \right)
\\
&=& \frac{q!}{\sigma^{2} N}\left( -2N\sum_{m = N }^{\infty}\left(\sum_{i = 1}^{\infty}\alpha_{i}\alpha_{i+\left|m\right|}\right)^{q} - 2\sum_{m =0 }^{N-1}m\left(\sum_{i = 1}^{\infty}\alpha_{i}\alpha_{i+\left|m\right|}\right)^{q}\right).
\end{eqnarray*}
By noticing that the condition $q > (2\beta -1)^{-1}$ is equivalent to $-q(2\beta -1) < -1$, we can apply Lemma \ref{l1} to get
\begin{eqnarray}
A_{q-1}(N) - 1 & \trianglelefteqslant & \sum_{m = N }^{\infty}m^{-q(2\beta -1)} + N^{-1}\sum_{m =0 }^{N-1}m^{-q(2\beta - 1) +1} \nonumber
\\
& \trianglelefteqslant & N^{-q(2\beta -1) + 1} + N^{-1}(1 + N^{-q(2\beta -1) + 2})\nonumber
\end{eqnarray}
and finally
\begin{eqnarray}
A_{q-1}(N) - 1 \trianglelefteqslant N^{-1} + N^{q-2q\beta+1}.  %N^{-q(2\beta -1)+1}
\end{eqnarray}
Let us now treat the second term of (\ref{toev}), i.e. $\sum_{r = 0}^{q-2}\mathbf{E}\left(A_{r}^{2}(N)\right)$. Here we can assume that $r \leq q-2$ is fixed. We have
\begin{eqnarray}
\mathbf{E}\left(A_{r}^{2}(N)\right) &=& \mathbf{E}\left(\frac{q^{2}r!^{2}}{\sigma^{4}N^{2}}\left(\begin{array}{cc}q-1 \\r \end{array}\right)^{4}\sum_{i,j,k,l = 1}^{N}\left\langle f_{k},f_{l} \right\rangle_{\mathcal{H}}^{r+1}\left\langle f_{i},f_{j} \right\rangle_{\mathcal{H}}^{r+1} \right. \nonumber
\\
&& \left. \times I_{2q-2 - 2r}\left(f_{k}^{\otimes q-1-r} \widetilde{\otimes}  f_{l}^{\otimes q-1-r}\right)I_{2q-2 - 2r}\left(f_{i}^{\otimes q-1-r} \widetilde{\otimes}  f_{j}^{\otimes q-1-r}\right)\right)\nonumber
\\
&=& c(r,q)N^{-2}\sum_{i,j,k,l = 1}^{N}\left\langle f_{k},f_{l} \right\rangle_{\mathcal{H}}^{r+1}\left\langle f_{i},f_{j} \right\rangle_{\mathcal{H}}^{r+1}\nonumber
\\
&& \times \left\langle f_{k}^{\otimes q-1-r} \widetilde{\otimes}  f_{l}^{\otimes q-1-r},f_{i}^{\otimes q-1-r} \widetilde{\otimes}  f_{j}^{\otimes q-1-r} \right\rangle_{\mathcal{H}^{\otimes 2q - 2r - 2}}\nonumber
\\
&=& \sum_{\underset{\alpha + \nu = q-r-1 }{\alpha, \nu \geq 0}}\sum_{\underset{\gamma + \delta = q-r-1}{\gamma, \delta \geq 0}}c(r,q,\alpha,\nu,\gamma,\delta)B_{r,\alpha,\nu,\gamma,\delta}(N)\nonumber
\end{eqnarray}
where
\begin{eqnarray*}
B_{r,\alpha,\nu,\gamma,\delta}(N) &=& N^{-2}\sum_{i,j,k,l = 1}^{N}\left\langle f_{k},f_{l} \right\rangle_{\mathcal{H}}^{r+1}\left\langle f_{i},f_{j} \right\rangle_{\mathcal{H}}^{r+1} \left\langle f_{k},f_{i} \right\rangle_{\mathcal{H}}^{\alpha}\left\langle f_{k},f_{j} \right\rangle_{\mathcal{H}}^{\nu}\left\langle f_{l},f_{i} \right\rangle_{\mathcal{H}}^{\gamma}\left\langle f_{l},f_{j} \right\rangle_{\mathcal{H}}^{\delta}\nonumber
\\
&=& N^{-2}\sum_{i,j,k,l = 1}^{N}\rho(k-l)^{r+1}\rho(i-j)^{r+1}\rho(k-i)^{\alpha}\rho(k-j)^{\nu}\rho(l-i)^{\gamma}\rho(l-j)^{\delta}.
\end{eqnarray*}
When $\alpha$,$\nu$,$\gamma$ and $\delta$ are fixed, we can decompose the sum $\sum_{i,j,k,l = 1}^{N}$ which appears in $B_{r,\alpha,\nu,\gamma,\delta}(N)$ just above, as follows:
\begin{eqnarray*}
\sum_{i=j=k=l} + \left( \sum_{\underset{l \neq i }{i=j=k}} + \sum_{\underset{k \neq i }{i=j=l}} + \sum_{\underset{j \neq i }{i=l=k}} + \sum_{\underset{i \neq j }{j=k=l}} \right) + \left( \sum_{\underset{k \neq i }{i=j,k=l}} + \sum_{\underset{j \neq i }{i=k,j=l}} + \sum_{\underset{j \neq i }{i=l,j=k}} \right)
\\
+ \left( \sum_{\underset{k\neq l,l\neq i}{i=j,k\neq i}} + \sum_{\underset{j\neq l,k\neq l}{i=k,j\neq i}} + \sum_{\underset{k\neq j,j\neq i}{i=l,k\neq i}} + \sum_{\underset{k\neq l,l\neq i}{j=k,k\neq i}}+ \sum_{\underset{k\neq l,l\neq i}{j=l,k\neq i}} + \sum_{\underset{k\neq j,j\neq i}{k=l,k\neq i}} \right) + \sum_{\underset{i \neq j \neq k \neq l}{i,j,k,l}}.
\end{eqnarray*}
We will have to evaluate each of these fifteen sums separatly. Before that, we will give a useful lemma that we will be using regularly throughout the paper.
\begin{lemma}
\label{lemmevaij}
For any $\alpha \in \mathbb{R}$, we have
\begin{eqnarray*}
\sum_{i \neq j = 1}^{n}\left|i-j\right|^{\alpha} = \sum_{i,j = 0}^{n-1}\left|i-j\right|^{\alpha} \trianglelefteqslant n\sum_{j = 0}^{n-1}j^{\alpha}.
\end{eqnarray*}
\end{lemma}
\begin{dem}
The following upperbounds prove this lemma
\begin{eqnarray*}
\left|\frac{\sum_{i,j = 0}^{n-1}\left|i-j\right|^{\alpha}}{n\sum_{j = 0}^{n-1}j^{\alpha}}\right| = \left|\frac{\sum_{m = 0}^{n-1}(n-m)m^{\alpha}}{n\sum_{j = 0}^{n-1}j^{\alpha}}\right| & \leq & \left|\frac{n\sum_{m = 0}^{n-1}m^{\alpha}}{n\sum_{j = 0}^{n-1}j^{\alpha}}\right| + \left|\frac{\sum_{m = 0}^{n-1}m^{\alpha+1}}{n\sum_{j = 0}^{n-1}j^{\alpha}}\right|
\\
& \leq & 1 + \left|\frac{\sum_{m = 0}^{n-1}m^{\alpha+1}}{\sum_{j = 0}^{n-1}j^{\alpha+1}}\right| \leq 2.
\end{eqnarray*}
\qed
\end{dem}
\\
\noindent Let's get back to our sums and begin by treating the first one. The first sum can rewritten as
\begin{eqnarray*}
&& N^{-2}\sum_{i=j=k=l}\rho(k-l)^{r+1}\rho(i-j)^{r+1}\rho(k-i)^{\alpha}\rho(k-j)^{\nu}\rho(l-i)^{\gamma}\rho(l-j)^{\delta}
\\
&&= N^{-2}\sum_{i=1}^{N}\rho(0)^{2r+2+\alpha+\nu+\gamma+\delta} = N^{-2}N \trianglelefteqslant  N^{-1}.
\end{eqnarray*}
For the second sum, we can write
\begin{eqnarray*}
&& N^{-2}\sum_{\underset{l \neq i }{i=j=k}}\rho(k-l)^{r+1}\rho(i-j)^{r+1}\rho(k-i)^{\alpha}\rho(k-j)^{\nu}\rho(l-i)^{\gamma}\rho(l-j)^{\delta}
\\
&&= N^{-2}\sum_{\underset{l \neq i }{i=j=k}}\rho(l-i)^{r+1+\gamma+\delta} = N^{-2}\sum_{i\neq l}\rho(l-i)^{q}.
\end{eqnarray*}
At this point, we will use lemma \ref{l1} and then lemma \ref{lemmevaij} to write
\begin{eqnarray*}
&& N^{-2}\sum_{\underset{l \neq i }{i=j=k}}\rho(k-l)^{r+1}\rho(i-j)^{r+1}\rho(k-i)^{\alpha}\rho(k-j)^{\nu}\rho(l-i)^{\gamma}\rho(l-j)^{\delta}
\\
&& \trianglelefteqslant N^{-2}\sum_{i\neq l = 1}^{N}\left|l-i\right|^{q(-2\beta +1)} \trianglelefteqslant N^{-1}\sum_{l = 1}^{N-1}l^{q(-2\beta +1)} \trianglelefteqslant N^{-1}(1 + N^{-2\beta q +q + 1})
\\
&& \trianglelefteqslant N^{-1} + N^{-2\beta q +q}.
\end{eqnarray*}
For the third sum, we are in the exact same case, therefore we obtain the same bound $N^{-1} + N^{-2\beta q +q}$. The fourth sum can be handled as follows
\begin{eqnarray*}
&& N^{-2}\sum_{\underset{j \neq i }{i=k=l}}\rho(k-l)^{r+1}\rho(i-j)^{r+1}\rho(k-i)^{\alpha}\rho(k-j)^{\nu}\rho(l-i)^{\gamma}\rho(l-j)^{\delta}
\\
&&= N^{-2}\sum_{\underset{j \neq i }{i=k=l}}\rho(i-j)^{r+1+\nu+\delta} \trianglelefteqslant N^{-2}\sum_{j\neq i}\left|i-j\right|^{(r+1+\nu+\delta)(-2\beta +1)}.
\end{eqnarray*}
Note that $r+1+\nu+\delta \geq 1$, so we get
\begin{eqnarray*}
&& N^{-2}\sum_{\underset{j \neq i }{i=k=l}}\rho(k-l)^{r+1}\rho(i-j)^{r+1}\rho(k-i)^{\alpha}\rho(k-j)^{\nu}\rho(l-i)^{\gamma}\rho(l-j)^{\delta}
\\
&&  \trianglelefteqslant N^{-2}\sum_{j\neq i}\left|i-j\right|^{-2\beta +1} \trianglelefteqslant N^{-1}\sum_{j=1}^{N-1}j^{-2\beta +1}
 \trianglelefteqslant N^{-1}(1+N^{-2\beta +2})
\\
&&  \trianglelefteqslant N^{-1}+N^{-2\beta +1}.
\end{eqnarray*}
For the fifth sum, we are in the exact same case and we obtain the same bound $N^{-1} + N^{-2\beta +1}$. For the sixth sum, we can proceed as follows
\begin{eqnarray*}
&& N^{-2}\sum_{\underset{k \neq i }{i=j,k=l}}\rho(k-l)^{r+1}\rho(i-j)^{r+1}\rho(k-i)^{\alpha}\rho(k-j)^{\nu}\rho(l-i)^{\gamma}\rho(l-j)^{\delta}
\\
&&= N^{-2}\sum_{k \neq i}\rho(k-i)^{\alpha+\nu+\gamma+\delta} = N^{-2}\sum_{k \neq i}\rho(k-i)^{2q-2r-2}.
\end{eqnarray*}
Recalling that $r \leq q-2 \Leftrightarrow 2(q-r-1) \geq 2$, we obtain
\begin{eqnarray*}
&& N^{-2}\sum_{\underset{k \neq i }{i=j,k=l}}\rho(k-l)^{r+1}\rho(i-j)^{r+1}\rho(k-i)^{\alpha}\rho(k-j)^{\nu}\rho(l-i)^{\gamma}\rho(l-j)^{\delta}
\\
&& \trianglelefteqslant N^{-2}\sum_{k \neq i}\left|k-i\right|^{(2q-2r-2)(-2\beta +1)}
\trianglelefteqslant N^{-2}\sum_{k \neq i}\left|k-i\right|^{-4\beta +2}
 \trianglelefteqslant N^{-1}\sum_{k = 1}^{N-1}k^{-4\beta +2}
\\
&& \trianglelefteqslant N^{-1} + N^{-4\beta +2}.
\end{eqnarray*}
We obtain the same bound, $N^{-1} + N^{-4\beta +2}$, for the seventh and eighth sums.  For the ninth sum, we have to deal with the following quantity.
\begin{eqnarray*}
&& N^{-2}\sum_{\underset{k\neq l,l\neq i }{i=j,k\neq i}}\rho(k-l)^{r+1}\rho(i-j)^{r+1}\rho(k-i)^{\alpha}\rho(k-j)^{\nu}\rho(l-i)^{\gamma}\rho(l-j)^{\delta}
\\
&&= N^{-2}\sum_{\underset{k\neq l,l\neq i }{k\neq i}}\rho(k-l)^{r+1}\rho(k-i)^{q-r-1}\rho(l-i)^{q-r-1}.
\end{eqnarray*}
For $\sum_{\underset{k\neq l,l\neq i }{k\neq i}}$, observe that it can be decomposed into
\begin{eqnarray}
\label{decompunesomme}
\sum_{k>l>i} + \sum_{k>i>l} + \sum_{l>i>k} + \sum_{i>l>k} + \sum_{i>k>l}.
\end{eqnarray}
For the first of the above sums, we can write
\begin{eqnarray*}
&& N^{-2}\sum_{k>l>i}\rho(k-l)^{r+1}\rho(k-i)^{q-r-1}\rho(l-i)^{q-r-1}
\\
&& \trianglelefteqslant  N^{-2}\sum_{k>l>i}(k-l)^{(r+1)(-2\beta +1)}(k-i)^{(q-r-1)(-2\beta +1)}(l-i)^{(q-r-1)(-2\beta +1)}
\\
&& \trianglelefteqslant  N^{-2}\sum_{k>l>i}(k-l)^{q(-2\beta +1)}(l-i)^{(q-r-1)(-2\beta +1)} \mbox{\  \  since\  } k-i > k-l
\\
&& =  N^{-2}\sum_{k}\sum_{l<k}(k-l)^{q(-2\beta +1)}\sum_{i<l}(l-i)^{(q-r-1)(-2\beta +1)}
\\
&& \trianglelefteqslant   N^{-2}\sum_{k}\sum_{l<k}(k-l)^{q(-2\beta +1)}\sum_{i<l}(l-i)^{-2\beta +1} \mbox{\  \  since\  } q-r-1 \geq 1
\\
&& \trianglelefteqslant   N^{-2}\sum_{k=1}^{N}\sum_{l=1}^{k-1}(k-l)^{q(-2\beta +1)}\sum_{i=1}^{l-1}(l-i)^{-2\beta +1}.
\end{eqnarray*}
Note that $\sum_{l=1}^{k-1}(k-l)^{q(-2\beta +1)} = \sum_{l=1}^{k-1}l^{q(-2\beta +1)}$ and that $\sum_{i=1}^{l-1}(l-i)^{-2\beta +1} = \sum_{i=1}^{l-1}i^{-2\beta +1}$. We can also bound the terms   $\sum_{l=1}^{k-1}l^{q(-2\beta +1)}$  (resp. $\sum_{i=1}^{l-1}i^{-2\beta +1}$) from above by $\sum_{l=1}^{N-1}l^{q(-2\beta +1)}$ (resp. $\sum_{i=1}^{N-1}i^{-2\beta +1}$). It follows that
\begin{eqnarray*}
&& N^{-2}\sum_{k>l>i}\rho(k-l)^{r+1}\rho(k-i)^{q-r-1}\rho(l-i)^{q-r-1}
\\
&& \trianglelefteqslant   N^{-2}\sum_{k=1}^{N}\sum_{l=1}^{N-1}l^{q(-2\beta +1)}\sum_{i=1}^{N-1}i^{-2\beta +1}
\\
&& \trianglelefteqslant   N^{-1}\sum_{l=1}^{N-1}l^{q(-2\beta +1)}\sum_{i=1}^{N-1}i^{-2\beta +1}
\\
&& \trianglelefteqslant   N^{-1}(1+N^{q(-2\beta +1)})(1+N^{-2\beta +1}) \trianglelefteqslant   N^{-1}(1+N^{-2\beta q +q +1})(1+N^{-2\beta +2})
\\
&& \trianglelefteqslant   N^{-1} + N^{-2\beta+1} + N^{-2\beta q + q} + N^{-2\beta q-2\beta +2}.
\end{eqnarray*}
Since $-2\beta+1 <0$, $-2\beta q + q <0$ and that $-2\beta q-2\beta +2 <0$, it is easy to check that $$-2\beta q-2\beta +2 < -2\beta q + q < -2\beta+1.$$ Consequently,
\begin{eqnarray*}
&& N^{-2}\sum_{k>l>i}\rho(k-l)^{r+1}\rho(k-i)^{q-r-1}\rho(l-i)^{q-r-1} \trianglelefteqslant  N^{-1} + N^{-2\beta+1}.
\end{eqnarray*}
We obtain the exact same bound $N^{-1} + N^{-2\beta+1}$ for the other terms of the decomposition (\ref{decompunesomme}) as well as for the tenth, eleventh, twelfth, thirteenth and fourteenth sums by applying the exact same method.
\\\\
\noindent This leaves us with the last (fifteenth) sum. We can decompose $\sum_{\underset{i \neq j \neq k \neq l}{i,j,k,l}}$ as follows
\begin{eqnarray}
\label{sommetoutediff}
\sum_{k>l>i>j} + \sum_{k>l>j>i}+ ...
\end{eqnarray}
For the first term, we have
\begin{eqnarray*}
&& N^{-2}\sum_{k>l>i>j}\rho(k-l)^{r+1}\rho(i-j)^{r+1}\rho(k-i)^{\alpha}\rho(k-j)^{\nu}\rho(l-i)^{\gamma}\rho(l-j)^{\delta}
\\
&& \trianglelefteqslant  N^{-2}\sum_{k>l>i>j}(k-l)^{q(-2\beta+1)}(i-j)^{(r+1)(-2\beta+1)}(l-i)^{(q-r-1)(-2\beta +1)}
\\
&& = N^{-2}\sum_{k}\sum_{l<k}(k-l)^{q(-2\beta+1)}\sum_{i<l}(l-i)^{(q-r-1)(-2\beta +1)}\sum_{j<i}(i-j)^{(r+1)(-2\beta+1)}
\\
&& \trianglelefteqslant  N^{-1}\sum_{l=1}^{N-1}l^{q(-2\beta+1)}\sum_{i=1}^{N-1}i^{(q-r-1)(-2\beta +1)}\sum_{j=1}^{N-1}j^{(r+1)(-2\beta+1)}
\\
&& \trianglelefteqslant  N^{-1}(1+N^{-2\beta q + q + 1})(1+N^{(q-r-1)(-2\beta +1) + 1})(1+N^{(r+1)(-2\beta+1) + 1})
\\
&& \trianglelefteqslant  N^{-1}(1 + N^{-2\beta q + q + 1})(1+N^{(r+1)(-2\beta + 1) + 1} + N^{q(-2\beta + 1) - (r+1)(-2\beta + 1) + 1} + N^{q(-2\beta + 1) + 2})
\\
&& \trianglelefteqslant  N^{-1}(1 + N^{-2\beta q + q + 1})(1+N^{-2\beta + 2} + N^{-2\beta + 2} + N^{q(-2\beta + 1) + 2}) \mbox{\  \  since\  } r+1,q-r-1 \geq 1
\\
&& \trianglelefteqslant  N^{-1}(1 + N^{-2\beta + 2} + N^{q(-2\beta + 1) + 2})
\\
&& \trianglelefteqslant  N^{-1} + N^{-2\beta + 1} + N^{q(-2\beta + 1) + 1}.
\end{eqnarray*}
We find the same bound $N^{-1} + N^{-2\beta + 1} + N^{q(-2\beta + 1) + 1}$ for the other terms of the decomposition (\ref{sommetoutediff}).
\\\\
\noindent Finally, by combining all these bounds, we find that
\begin{eqnarray*}
\underset{r=1,...,q-1}{\mbox{max}} \mathbf{E}\left(A_{r}^{2}\right) \trianglelefteqslant N^{-2\beta + 1} + N^{q(-2\beta + 1) + 1},
\end{eqnarray*}
and we obtain
\begin{eqnarray*}
\mathbf{E}\left(\left(\frac{1}{q}\left\|DZ_{N}\right\|_{\mathcal{H}}^{2} - 1\right)^{2}\right) \trianglelefteqslant N^{-2\beta + 1} + N^{q(-2\beta + 1) + 1},
\end{eqnarray*}
which allow us to complete the proof.
\qed
\begin{remark}
\begin{enumerate}
   \item When $q = 2$, $\frac{q}{2q-2} = \frac{1}{2}$, so the first line vanishes. If $q > 2$, both lines exists and $\frac{q}{2q-2} \underset{q \rightarrow
+\infty}{\longrightarrow} \frac{1}{2}$.
   \item When $q < (2\beta -1)^{-1}$, the sequence $Z_{N}$ does not converge in law towards an $\mathcal{N}(0,1)$. It converges (with another normalization) to a Hermite random variable.

       \item The results in the above theorem are coherent with those found in \cite{NoPe1}, Theorem 4.1. Indeed, in \cite{NoPe1} one works with $Y_{n}=B^{H}_{n+1}-B^{H}_{n}$ instead of $X_{n}$, where $B^{H}$ is a fractional Brownian motion. Note that the covariance function $\rho'(m)=\mathbf{E}\left(Y_{0}Y_{m}\right)$ of $Y$ behaves as $m^{2H-2}$ while, as it follows from Lemma \ref{l1}, the covariance  of $X$ behaves as $m^{-2\beta+1}$. Thus $\beta$ corresponds to $\frac{3}{2}-H$. It can be seen that Theorem \ref{tt1} is in concordance with Theorem 4.1 in \cite{NoPe1}.
\end{enumerate}
\end{remark}

\subsection{Error Bounds in the Non-Central Limit Theorem}

We will now turn our attention to the case where $q < (2\beta -1)^{-1}$, where we will use the total variation distance instead of the Kolmogorov distance because that is the distance which appears in a result by Davydov and Martynova \cite{DaMa}. This result will be central to our proof of the bounds. Recall that the total variation distance between two real-valued random variables X and Y probability distributions is defined by
\begin{eqnarray}
d_{\mbox{\tiny{TV}}}(\mathcal{L}(X),\mathcal{L}(Y)) = \underset{A \in \mathcal{B}(\mathbb{R})}{\mbox{sup}} \left|\mathbf{P}(Y \in A) - \mathbf{P}(X \in A)\right|
\end{eqnarray}
where $\mathcal{B}(\mathbb{R})$ denotes the class of Borel sets of $\mathbb{R}$. We have the following result by Davydov and Martynova \cite{DaMa} on the total variation distance between elements of a fixed Wiener chaos.
\begin{theorem}
\label{davmart}
Fix an integer $q \geq 2$ and let $f \in \mathcal{H}^{\odot q} \backslash \left\{0\right\}$. Then, for any sequence $\left\{f_{n}\right\}_{n \geq 1} \subset \mathcal{H}^{\odot q}$ converging to $f$, their exists a constant $c_{q,f}$, depending only on $q$ and $f$, such that
\begin{eqnarray*}
d_{\mbox{\tiny{TV}}}(I_{q}(f_{n}),I_{q}(f)) \leq c_{q,f} \left\|f_{n} - f\right\|_{\mathcal{H}^{\odot q}}^{1/q}.
\end{eqnarray*}
\end{theorem}
~\\
\noindent We will now use the scaling property of the Brownian motion to introduce a new sequence $U_{N}$ that has the same law as $S_{N}$. Recall that $S_{N}$ is defined by
\begin{eqnarray*}
S_{N} = \sum_{n=1}^{N}H_{q}\left(\sum_{i=1}^{\infty}\alpha_{i}\left(W_{n-i-1} - W_{n-i}\right)\right).
\end{eqnarray*}
Let $U_{N}$ be defined by
\begin{eqnarray*}
U_{N} = \sum_{n=1}^{N}H_{q}\left(\sum_{i=1}^{\infty}\alpha_{i}N^{\frac{1}{2}}\left(W_{\frac{n-i-1}{N}} - W_{\frac{n-i}{N}}\right)\right).
\end{eqnarray*}
Based on the scaling property of the Brownian motion, $U_{N}$ has the same law as $S_{N}$ for every fixed $N$. Recall that Theorem \ref{CLTGen} states that
\begin{equation*}
h_{q, \beta}^{-1}N^{\beta q-\frac{q}{2}-1}S_{N} \underset{N \rightarrow +\infty}{\longrightarrow} Z^{(q)}
\end{equation*}
where $Z^{(q)}$ is a Hermite random variable of order $q$ (it is actually the value at time 1 of the Hermite process of order $q $ with self-similarity index
$$\frac{q}{2}-q\beta +1$$
defined in \cite{CNT}).
Let us first prove the following renormalization result.

\begin{lemma}
\label{reno}
Let
\begin{equation}
\label{hq}
h_{q, \beta}^{2}= \frac{2c_{\beta}^{q}}{ q! (-2\beta q +q+1)(-2\beta +q+2)}.
\end{equation}Then
\begin{equation*}
\mathbf{E} \left( h_{q, \beta}^{-1} N ^{\beta q -\frac{q}{2}-1} S_{N} \right) ^{2} \underset{N \rightarrow +\infty}{\longrightarrow} 1.
\end{equation*}
\end{lemma}
{\bf Proof: }
Indeed, since $S_{N}=\frac{1}{q!}I_{q}(f_{N})$ we have
\begin{eqnarray*}
\mathbf{E} \left( h_{q, \beta}^{-1} N ^{\beta q -\frac{q}{2}-1} S_{N} \right) ^{2}&=&h_{q, \beta}^{-2}  \frac{1}{(q!)}N^{2\beta q-q-2}\sum_{n,m=1}^{N}\rho(\vert n-m\vert )^{q}\\
&=&h_{q, \beta}^{-2}  \frac{1}{(q!)}N^{2\beta q-q-2}N\rho(0)^{q} + 2h_{q, \beta}^{-2}  \frac{1}{(q!)}N^{2\beta q-q-2}\sum_{n,m=1; n>m}^{N}\rho( n-m )^{q} \\
&\sim &2h_{q, \beta}^{-2}  \frac{1}{(q!)^{2}}N^{2\beta q-q-2}\sum_{n,m=1; n>m}^{N}\rho( n-m )^{q}
\end{eqnarray*}
where for the last equivalence we notice that the diagonal term $h_{q, \beta}^{-2}  \frac{1}{(q!)}N^{2\beta q-q-2}N\rho(0)^{q}$ converges to zero since $q<\frac{1}{2\beta -1}$. Therefore, by using the change of indices $n-m=k$ we can write
\begin{eqnarray*}
\mathbf{E} \left( h_{q, \beta}^{-1} N ^{\beta q -\frac{q}{2}-1} S_{N} \right) ^{2}&=&h_{q, \beta}^{-2}  \frac{1}{(q!)}N^{2\beta q-q-2}\sum_{n,m=1}^{N}\rho(\vert n-m\vert )^{q}\\
&\sim &2h_{q, \beta}^{-2}  \frac{1}{(q!)}N^{2\beta q-q-2}\sum_{k=1}^{N} (N-k) \rho (k)^{q} \\
&\sim & 2h_{q, \beta}^{-2}  \frac{1}{(q!)}N^{2\beta q-q-2}\sum_{k=1}^{N} (N-k)k^{-2\beta q+q}
\end{eqnarray*}
because,  according to Lemma \ref{l1}, $\rho(k)$ behaves as  $c_{\beta} k^{-2\beta +1} $ when $k$ goes to $\infty$. Consequently,
\begin{equation*}
\mathbf{E} \left( h_{q, \beta}^{-1} N ^{\beta q -\frac{q}{2}-1} S_{N} \right) ^{2}\sim  2h_{q, \beta}^{-2}  \frac{1}{q!}\frac{1}{N}\sum_{k=1}^{N} \left( 1-\frac{k}{N}\right) \left( \frac{k}{N}\right) ^{-2\beta q+q}
\end{equation*}
and this converges to 1 as $N\to \infty$  because  $\frac{1}{N}\sum_{k=1}^{N} \left( 1-\frac{k}{N}\right) \left( \frac{k}{N}\right) ^{-2\beta q+q} $ converges to
$$\int_{0}^{1} (1-x) x^{-2\beta q+q} dx = \frac{1}{(-2\beta q+q+1)(-2\beta q+q+2)}.$$\qed
\\\\
\noindent Let $Z_{N}$ be defined here by
\begin{eqnarray*}
Z_{N} = N^{\beta q - \frac{q}{2} - 1} U_{N} = N^{\beta q - \frac{q}{2} - 1}\sum_{n=1}^{N}H_{q}\left(\sum_{i=1}^{\infty}\alpha_{i}N^{\frac{1}{2}}\left(W_{\frac{n-i-1}{N}} - W_{\frac{n-i}{N}}\right)\right).
\end{eqnarray*}
We also know that $h_{q,\beta}^{-1} Z_{N} \underset{N \rightarrow +\infty}{\longrightarrow} Z^{(q)} $  in law (because $U_{N}$ has the same law as $S_{N}$), with $Z^{(q)}$ given by (\ref{hermite}). Let us give a proper representation of $Z_{N}$ as an element of the $q^{\mbox{\tiny{th}}}$-chaos. We have
\begin{eqnarray*}
Z_{N} & = & N^{\beta q - \frac{q}{2} - 1}\sum_{n=1}^{N}H_{q}\left(\sum_{i=1}^{\infty}\alpha_{i}N^{\frac{1}{2}}\left(W_{\frac{n-i-1}{N}} - W_{\frac{n-i}{N}}\right)\right)
\\
& = & N^{\beta q - \frac{q}{2} - 1}\sum_{n=1}^{N}H_{q}\left(I_{1}\left(N^{\frac{1}{2}}\sum_{i=1}^{\infty}\alpha_{i}\mathbf{1}_{\left[\frac{n-i-1}{N}, \frac{n-i}{N}\right]}\right)\right)
\\
& = & N^{\beta q - \frac{q}{2} - 1}\sum_{n=1}^{N}\frac{1}{q!}I_{q}\left(\left(N^{\frac{1}{2}}\sum_{i=1}^{\infty}\alpha_{i}\mathbf{1}_{\left[\frac{n-i-1}{N}, \frac{n-i}{N}\right]}\right)^{\otimes q}\right)
\\
& = &\frac{1}{q!} I_{q}\left(N^{\beta q - 1}\sum_{n=1}^{N}\left(\sum_{i=1}^{\infty}\alpha_{i}\mathbf{1}_{\left[\frac{n-i-1}{N}, \frac{n-i}{N}\right]}\right)^{\otimes q}\right)
\\
& := &\frac{1}{q!} I_{q}\left(\underbrace{N^{\beta q - 1}\sum_{n=1}^{N}g_{n}^{\otimes q}}_{g_{N}}\right)
\end{eqnarray*}
with $g_{n} = \sum_{i=1}^{\infty}\alpha_{i}\mathbf{1}_{\left[\frac{n-i-1}{N}, \frac{n-i}{N}\right]}$ and $g_{N} = N^{\beta q - 1}\sum_{n=1}^{N}g_{n}^{\otimes q}  \in \mathcal{H}^{\odot q}$. We will see that $h_{q,\beta }^{-1}Z_{N}$ converges towards $Z$ in $L^{2}(\Omega)$, or equivalently that $\left\{ \frac{1}{q!}h_{q,\beta }^{-1}g_{N}\right\}_{N \geq 1}$  converges in $L^{2}(\mathbb{R}^{\otimes q})={\cal{H}}^{\otimes q}$ to the kernel $g$ of the Hermite random variable (\ref{kernel}) by computing the following $L^{2}$ norm.
\begin{eqnarray*}
\mathbf{E}\left(\left|h_{q,\beta}^{-1}Z_{N} - Z\right|^{2}\right) = \mathbf{E}\left(\left|I_{q}(\frac{1}{q!}h_{q,\beta }^{-1}g_{N}) - I_{q}(g)\right|^{2}\right) = q!\left\|\frac{1}{q!}h_{q,\beta }^{-1}g_{N} - g\right\|_{\mathcal{H}^{\otimes q}}^{2}.
\end{eqnarray*}
We will now study $\left\|g_{N} - g\right\|_{\mathcal{H}^{\otimes q}}^{2}$ and establish the rate of convergence of this quantity.

\begin{prop}\label{pp1}
We have
\begin{equation*}
\left\|h_{q, \beta}^{-1}\frac{1}{q!}g_{N} - g\right\|_{\mathcal{H}^{\otimes q}}^{2}={\cal{O}}(N^{2\beta q-q-1}).
\end{equation*}
In particular the sequence $h_{q, \beta}^{-1}\frac{1}{q!}g_{N}$ converges in $L^{2}(\mathbb{R}^{\otimes q})$ as $N\to \infty$ to the kernel of the Hermite process $g$ (\ref{kernel}).
\end{prop}
{\bf Proof: } We have
\begin{eqnarray}
\left\|g_{N} \right\|_{\mathcal{H}^{\otimes q}}^{2} & = & N^{2\beta q - 2}\sum_{n,k=1}^{N}\left\langle g_{n},g_{k}\right\rangle_{\mathcal{H}}^{q}
\nonumber\\
& = & N^{2\beta q - 2}\sum_{n,k=1}^{N}\left(\int_{\mathbb{R}}\sum_{i=1}^{\infty}\sum_{j=1}^{\infty}\alpha_{i}\alpha_{j}\mathbf{1}_{\left[\frac{n-i-1}{N}, \frac{n-i}{N}\right]}(u)\mathbf{1}_{\left[\frac{k-j-1}{N}, \frac{k-j}{N}\right]}(u)du\right)^{q}
\nonumber\\
& = & N^{2\beta q - 2}\sum_{n,k=1}^{N}\left(\sum_{i=1}^{\infty}\alpha_{i}\alpha_{i + \left|n-k\right|}\int_{\frac{n-i-1}{N}}^{\frac{n-i}{N}}du\right)^{q}
\nonumber\\
& = & N^{2\beta q - q - 2}\sum_{n,k=1}^{N}\left(\sum_{i=1}^{\infty}i^{-\beta}\left(i + \left|n-k\right|\right)^{-\beta}\right)^{q}.\label{estim1}
\end{eqnarray}
\\\\
\noindent In addition, based on the definition of the Hermite process, we have
$$q!\left\|g \right\|_{\mathcal{H}^{\otimes q}}^{2} =1.$$
\noindent Let us now compute the scalar product $\langle g_{N}, g\rangle _{{\cal{H}}^{\otimes q}}$ where $g$ is given by (\ref{kernel}). It holds that
\begin{eqnarray*}
&&\langle g_{N}, g\rangle _{{\cal{H}}^{\otimes q}}=d(q,\beta) N^{\beta q-1} \sum_{n=1}^{N} \langle g_{n} ^{\otimes q} , g\rangle _{{\cal{H}}^{\otimes q}}\\
&=& d(q,\beta) N^{\beta q-1} \sum_{n=1}^{N}\int_{0}^{1}  \left( \sum_{i\geq 1} \alpha _{i} \int_{\mathbb{R}} (u-y)_{+}^{-\beta } 1_{\left( \frac{n-i-1}{N}, \frac{n-i}{N}\right] }(y)dy \right) ^{q}du\\
&=&d(q,\beta) N^{\beta q-1} \sum_{n=1}^{N}\sum_{k=1}^{N} \int_{\frac{k-1}{N}}^{\frac{k}{N}} \left( \sum_{i\geq 1} \alpha _{i} \int_{\mathbb{R}} (u-y)_{+}^{-\beta } 1_{\left( \frac{n-i-1}{N}, \frac{n-i}{N}\right] }(y)dy \right) ^{q}du.
\end{eqnarray*}
We will now perform the change of variables $u'=(u-\frac{k-1}{N})N$ and $y'=(y-\frac{n-i-1}{N})N$ (renaming the variables by $u$ and $y$), obtaining
\begin{eqnarray*}
\langle g_{N}, g\rangle _{{\cal{H}}^{\otimes q}}&=&d(q,\beta) N^{\beta q-1} N^{-q-1}\sum_{n=1}^{N}\sum_{k=1}^{N} \int_{0}^{1} \left( \sum_{i\geq 1} \alpha _{i}\int_{0}^{1} \left( \frac{u-y+k-n+i}{N}\right) _{+}^{-\beta } dy \right) ^{q}du \\
&\sim & d(q,\beta) N^{\beta q-q-2}\sum_{n=1}^{N}\sum_{k=1}^{N-1}\left( \sum_{i\geq 1} \alpha _{i}\left( \frac{k-n+i}{N}\right) _{+}^{-\beta } \right) ^{q}
\end{eqnarray*}
where we used the fact that when $N\to \infty$, the quantity $\frac{u-y}{N} $ is negligible. Hence, by eliminating the  diagonal term as above,
\begin{eqnarray*}
\langle g_{N}, g\rangle _{{\cal{H}}^{\otimes q}}&\sim& d(q,\beta) N^{2\beta q-q-2}\sum_{k,n=1; k>n} \left( \sum_{i\geq 1} \alpha _{i}  (i+k-n)^{-\beta } \right) ^{q}\\
&&+d(q,\beta) N^{2\beta q-q-2}\sum_{k,n=1; k<n} \left( \sum_{i\geq n-k} \alpha _{i}  (i+k-n)^{-\beta } \right) ^{q}
\end{eqnarray*}
and by using the change of indices $k-n=l$ in the first summand above and $n-k=l$ in the second summand we observe that
\begin{eqnarray}
\langle g_{N}, g\rangle _{{\cal{H}}^{\otimes q}}&\sim &d(q,\beta) N^{2\beta q-q-2}\sum_{l=1}^{N}(N-l) \left( \sum_{i\geq 1} i^{-\beta } (i+l) ^{-\beta }  \right) ^{q}\nonumber\\
&&+ d(q,\beta) N^{2\beta q-q-2}\sum_{l=1}^{N}(N-l) \left( \sum_{i\geq l} i^{-\beta } (i-l) ^{-\beta }  \right) ^{q}.\label{estim3}
\end{eqnarray}
By summarizing the above estimates (\ref{estim1}) and (\ref{estim3}), we establish that
\begin{eqnarray*}
\left\|h_{q, \beta}^{-1}\frac{1}{q!}g_{N} - g\right\|_{\mathcal{H}^{\otimes q}}^{2}&\sim &N ^{2\beta q -q-1} \left[ 2h_{q, \beta}^{-2}\frac{1}{(q!)^{2}} \frac{1}{N} \sum_{k=1}^{N} (N-k) \left( \sum_{i\geq 1} i^{-\beta } (i+k) ^{-\beta }  \right) ^{q} \right. \\
&&\left.  -2d(q, \beta) h_{q, \beta}^{-1}  \frac{1}{q!} \frac{1}{N} \sum_{k=1}^{N} (N-k) \left( \sum_{i\geq 1} i^{-\beta } (i+k) ^{-\beta }  \right) ^{q} \right. \\
&&\left. -2d(q, \beta) h_{q, \beta}^{-1}   \frac{1}{N} \sum_{k=1}^{N} (N-k) \left( \sum_{i\geq k} i^{-\beta } (i-k) ^{-\beta }  \right) ^{q}+ \frac{1}{q!}N^{-2\beta q +q+1} \right].
\end{eqnarray*}
To obtain the conclusion, it suffices to check that the sequence
\begin{eqnarray*}
a_{N}&:=&2h_{q, \beta}^{-2}\frac{1}{(q!)^{2}} \frac{1}{N} \sum_{k=1}^{N} (N-k) \left( \sum_{i\geq 1} i^{-\beta } (i+k) ^{-\beta }  \right) ^{q}\\
&&
  -2d(q, \beta) h_{q, \beta}^{-1} \frac{1}{q!}  \frac{1}{N} \sum_{k=1}^{N} (N-k) \left( \sum_{i\geq 1} i^{-\beta } (i+k) ^{-\beta }  \right) ^{q} \\
&& -2d(q, \beta) h_{q, \beta}^{-1}  \frac{1}{q!} \frac{1}{N} \sum_{k=1}^{N} (N-k) \left( \sum_{i\geq k} i^{-\beta } (i-k) ^{-\beta }  \right) ^{q}+ \frac{1}{q!} N^{-2\beta q +q+1}
\end{eqnarray*}
is uniformly bounded by a constant with respect to $N$. Since $d(q,\beta)h_{q, \beta}^{-1}= \frac{1}{q!} h_{q, \beta}^{-2}$, $\sum_{i\geq 1} i^{-\beta } (i+k) ^{-\beta }\sim c_{\beta }k^{-2\beta q+q }$  and
$$ \sum_{i\geq k} i^{-\beta } (i-k) ^{-\beta }= \sum_{i\geq 1} i^{-\beta } (i+k) ^{-\beta }$$
(by the change of notation $i-k=j$), the sequence $a_{N}$ can be written as
\begin{eqnarray*}
a_{N}\sim \frac{1}{q!}\left( - (-2\beta q+q+1)(-2\beta q +q+2 )\frac{1}{N} \sum_{k=1}^{N} (N-k)k^{-2\beta q+q}+ N^{-2\beta q +q+1}\right).
\end{eqnarray*}
It is easy to check that
\begin{eqnarray*}
N^{-2\beta q +q+1}&=& N^{-2\beta q +q+1} (-2\beta q+q+1)(-2\beta q +q+2 )\int_{0}^{1}(1-x) x^{-2\beta q +q}dx \\&=& (-2\beta q+q+1)(-2\beta q +q+2 )\frac{1}{N} \int_{0}^{N} (N-y) y^{-2\beta q+q }dy
\end{eqnarray*}
(by the change of variables $xN=y$). Thus,
\begin{eqnarray*}
q! a_{N}&\sim &c \frac{1}{N} \sum_{k=1}^{N} \int_{k-1} ^{k} dy\left( (N-y) y^{-2\beta q+q }-(N-k) k^{-2\beta q+q }\right)\\
& \leq & \sum_{k=1}^{N} \int_{k-1} ^{k} dy\left| y^{-2\beta q +q} -k^{-2\beta q+q} \right| +\frac{1}{N} \sum_{k=1}^{N} \int_{k-1} ^{k} dy\left| y^{-2\beta q +q+1} -k^{-2\beta q+q+1} \right| \\
& \leq &\sum_{k=1}^{N} \left( (k-1) ^{-2\beta q+q} -k^{-2\beta q+q} \right) + \frac{1}{N} \sum_{k=1}^{N} \left( k ^{-2\beta q+q+1} -(k-1)^{-2\beta q+q+1} \right)
\end{eqnarray*}
and elementary computations show that the terms on the last line above  are of order of $N^{-2\beta q+q+1}$. \qed
\\\\
\noindent As a consequence of Proposition \ref{pp1} and of Theorem \ref{davmart}, we obtain
\begin{theorem}
\label{tt2}
Let $q<\frac{1}{2\beta -1}$ and let $S_{N}$ be given by (\ref{sn}).
$$d_{TV}\left( h_{q,\beta }^{-1}N^{\beta q-\frac{q}{2}-1}S_{N}, Z^{(q)} \right)\leq C_{0}(q, \beta)N^{2\beta q -q-1}$$
where $Z^{(q)}$ is given by (\ref{hermite}), $h_{q,\beta}$ is given by (\ref{hq}) and $C_{0}(q, \beta)$ is a positive constant.
\end{theorem}

\section{Application: Hsu-Robbins and Spitzer  theorems for moving averages}
In this section, we will give an application of the bounds obtained in Theorems \ref{tt1} and \ref{tt2}. The purpose of the Spitzer theorem for moving averages is to find the asymptotic behavior as $\varepsilon \to 0$ of  the sequences
\begin{equation*}
f_{1}(\varepsilon) = \sum_{N\geq 1} \frac{1}{N} P\left( \left|S_{N}\right|> \varepsilon N\right).
\end{equation*}
when $q>\frac{1}{2\beta -1}$ and
\begin{equation*}
f_{2}(\varepsilon) = \sum_{N\geq 1} \frac{1}{N} P\left( \left|S_{N}\right|> \varepsilon N^{-2\beta q + q +2}\right).
\end{equation*}
when $q<\frac{1}{2\beta -1}$.
The cases of the increments of the fractional Brownian motion were treated in \cite{T1}. The same arguments can be applied here. Let us briefly describe the method used to find the limit of $f(\varepsilon)$ as $\eps \to 0$. Let $q>\frac{1}{2\beta -1}$ so the limit of $\sigma ^{-1}\frac{1}{\sqrt{N} } S_{N}$ is a standard normal random variable. We have
\begin{eqnarray*}
f_{1}(\eps)&=& \sum_{N\geq 1} \frac{1}{N} P\left(\sigma ^{-1}\frac{1}{\sqrt{N} } \left|S_{N}\right|> \frac{\eps  \sqrt{N}}{\sigma} \right)\\
&=& \sum_{N\geq 1} \frac{1}{N} P\left(\left|Z\right| > \frac{\eps  \sqrt{N}}{\sigma }  \right)\\
&&+ \sum_{N\geq 1} \frac{1}{N} \left[ P\left( \sigma^{-1} \frac{1}{\sqrt{N} } \left|S_{N}\right|> \frac{\eps  \sqrt{N}}{\sigma} \right)-P\left(\left|Z\right| > \frac{\eps  \sqrt{N}}{\sigma } \right)\right]
\end{eqnarray*}
where $Z$ denotes a standard normal random variable. The first summand above was estimated in \cite{T1},  Lemma 1 while the second summand converges to zero by using the bound in Theorem \ref{tt1} and the proof of the Proposition 1 in \cite{T1}. When $q<\frac{1}{2\beta -1}$, similarly
\begin{eqnarray*}
f_{2}(\eps)&=& \sum_{N\geq 1} \frac{1}{N}P\left( h_{q, \beta}^{-1} N^{\beta q -\frac{q}{2}-1}\left|S_{N}\right| > d_{q,\beta }^{-1} \eps N^{-\beta q+\frac{q}{2}+1} \right)\\
&=&\sum_{N\geq 1} \frac{1}{N}P\left( \left|Z^{(q)}\right|>d_{q,\beta }^{-1} \eps N^{-\beta q+\frac{q}{2}+1} \right)\\
&&+ \sum_{N\geq 1} \frac{1}{N}\left[ P\left( d_{q, \beta}^{-1} N^{\beta q -\frac{q}{2}-1}\left|S_{N}\right| > d_{q,\beta }^{-1} \eps N^{-\beta q+\frac{q}{2}+1} \right)-P\left( \left|Z^{(q)}\right|>d_{q,\beta }^{-1} \eps N^{-\beta q+\frac{q}{2}+1} \right) \right]
\end{eqnarray*}
with $Z^{(q)}$ a Hermite random variable of order $q$. The first summand was also estimated in \cite{T1}, Lemma 1 while the second summand can be handled as in Proposition 2 in \cite{T1} and the result in Theorem \ref{tt2}. Hence, we obtain
\begin{prop}
When $q>\frac{1}{2\beta -1}$,  $$\lim\limits_{\eps \to 0}\frac{1}{-\log (\eps) } f_{1}(\eps)=2$$ and when  $q<\frac{1}{2\beta -1}$  then $$\lim\limits_{\eps \to 0}\frac{1}{-\log (\eps) } f_{2}(\eps)=\frac{1}{1+\frac{q}{2}-\beta q }.$$
\end{prop}
It is also possible to give Hsu-Robbins type results, meaning to find the asymptotic behavior as $\eps \to 0$ of
\begin{equation*}
g_{1}(\varepsilon) = \sum_{N\geq 1}  P\left( \left|S_{N}\right|> \varepsilon N\right)
\end{equation*}
when $q>\frac{1}{2\beta -1}$ and
\begin{equation*}
g_{2}(\varepsilon) = \sum_{N\geq 1}  P\left( \left|S_{N}\right|> \varepsilon N^{-2\beta q + q +2}\right)
\end{equation*}
when $q<\frac{1}{2\beta -1}$.
This also follows from Section 4 in \cite{T1} and Theorems \ref{tt1} and \ref{tt2}.
\begin{prop}
When $q>\frac{1}{2\beta -1}$,  $$\lim\limits_{\eps \to 0} (\sigma_{q,\beta}^{-1} \eps)^{2}g_{1}(\eps)=1=\mathbf{E}\left(Z^{2}\right)$$ and when  $q<\frac{1}{2\beta -1}$  then $$\lim\limits_{\eps \to 0}(h_{q\beta}^{-1} \eps)^{\frac{1}{1+\frac{q}{2} - \beta q}}g_{2}(\eps)=\mathbf{E}\left| Z^{(q)}\right|^{\frac{1}{1+\frac{q}{2} - \beta q}}.$$
\end{prop}

\end{document}